\documentclass[12pt]{article}
\usepackage{amscd}
\begin{document}
\renewcommand{\thesubsection}{\arabic{subsection}}
\newenvironment{proof}{{\bf Proof}:}{\vskip 5mm }
\newenvironment{rem}{{\bf Remark}:}{\vskip 5mm }
\newenvironment{remarks}{{\bf Remarks}:\begin{enumerate}}{\end{enumerate}}
\newenvironment{examples}{{\bf Examples}:\begin{enumerate}}{\end{enumerate}}  
\newtheorem{proposition}{Proposition}[subsection]
\newtheorem{lemma}[proposition]{Lemma}
\newtheorem{definition}[proposition]{Definition}
\newtheorem{theorem}[proposition]{Theorem}
\newtheorem{cor}[proposition]{Corollary}
\newtheorem{conjecture}{Conjecture}
\newtheorem{pretheorem}[proposition]{Pretheorem}
\newtheorem{hypothesis}[proposition]{Hypothesis}
\newtheorem{example}[proposition]{Example}
\newtheorem{remark}[proposition]{Remark}
\newcommand{\llabel}[1]{\label{#1}}
\newcommand{\comment}[1]{}
\newcommand{\sr}{\rightarrow}
\newcommand{\dw}{\downarrow}
\newcommand{\bdl}{\bar{\Delta}}
\newcommand{\zz}{{\bf Z\rm}}
\newcommand{\zq}{{\bf Z}_{qfh}}
\newcommand{\nn}{{\bf N\rm}}
\newcommand{\qq}{{\bf Q\rm}}
\newcommand{\nq}{{\bf N}_{qfh}}
\newcommand{\oo}{\otimes}
\newcommand{\uu}{\underline}
\newcommand{\ih}{\uu{Hom}}
\newcommand{\af}{{\bf A}^1}
\newcommand{\dsr}{\stackrel{\sr}{\scriptstyle\sr}}
\begin{center}
{\Large\bf Unstable motivic homotopy categories in Nisnevich and
cdh-topologies.}\\
\vskip 4mm
\vskip 3mm

 {\bf Vladimir Voevodsky}\footnote{Supported by the NSF
grants DMS-97-29992 and DMS-9901219, Sloan Research Fellowship and Veblen
Fund}$^,$\footnote{School of Mathematics, Institute for Advanced Study,
Princeton NJ, USA. e-mail: vladimir@ias.edu}\\
{August 2000}
\end{center}
\vskip 4mm
\tableofcontents
\subsection{Introduction}
One can do the motivic homotopy theory in the context of different
motivic homotopy categories. One can vary the topology on the category
of schemes used to define the homotopy category or one can vary the
category of schemes itself considering only schemes satisfying certain
conditions. The category obtained by taking smooth schemes and the
Nisnevich topology seems to play a distinguished role in the theory
because of the Gluing Theorem (see \cite{MoVo}) and some other, less
significant, nice properties. On the other hand, in the parts of the
motivic homotopy theory dealing with the motivic cohomology it is
often desirable to be able to work with all schemes instead of just
the smooth ones. For example, the motivic Eilenberg-MacLane spaces are
naturally representable (in characteristic zero) by singular schemes
built out of symmetric products of projective spaces but we do not
know of any explicit way to represent these spaces by simplicial
smooth schemes.

The goal of this paper is to show that, under the resolution of
singularities assumption, the pointed motivic homotopy category of
smooth schemes over a field with respect to the Nisnevich topology is
almost equivalent to the pointed motivic homotopy category of all
schemes over the same field with respect to the cdh-topology. More
precisely, we show that the inverse image functor
$${\bf L}\pi^*:H_{\bullet}((Sm/k)_{Nis},\af)\sr
H_{\bullet}((Sch/k)_{cdh},\af)$$
from the former category to the later one is a localization and if $f$
is a morphism such that ${\bf L}\pi^*(f)$ is an isomorphism then the
first simplicial suspension of $f$ is an isomorphism. This should
imply in particular that the corresponding s-stable and T-stable
motivic homotopy categories are equivalent.

The present paper is a continuation of the series started with
\cite{HH0} and \cite{HH1} and it uses the formalism developed
there. In the first section we define the standard cd-structures on
the category of Noetherian schemes and prove that they are complete,
regular and bounded. In the next section we prove some simple results
about the homotopy categories of sites with interval with completely
decomposable topologies and apply them to get an explicit description
of the $\af$-weak equivalences in terms of $\Delta$-closed classes.
Our results also imply that the motivic homotopy categories defined
with respect to the standard topologies are homotopy categories of
almost finitely generated closed model structures (see \cite{HH0}). In
the last section we apply these results to prove the comparison
theorem.

This paper was written while I was a member of the Institute for
Advanced Study in Princeton and, part of the time, an employee of the
Clay Mathematics Institute. I am very grateful to both institutions
for their support. I would also like to thank Charles Weibel who
pointed out a number of places in the previous version of the paper
which required corrections.

\vskip 3mm

\noindent
Everywhere below a scheme means a Noetherian scheme.

\subsection{The standard cd-structures on categories of schemes}
\label{standard}
Let us consider the following two cd-structures on the category of
Noetherian schemes.
\begin{description}
\item[Upper cd-structure] or Nisnevich cd-structure where a square of the
form 
\begin{equation}
\llabel{eq1old}\llabel{eq1}
\begin{CD}
B@>>>Y\\
@VVV @VVpV\\
A @>e>> X
\end{CD}
\end{equation}
is distinguished if it is a pull-back square such
that $p$ is etale, $e$ is an open embedding and $p^{-1}(X-e(A))\sr
X-e(A)$ is an isomorphism. Here $X-e(A)$ is considered with the
reduced scheme structure. 
\item[Lower cd-structure] or proper cdh-structure where a square of the
form (\ref{eq1old}) is distinguished if it is a pull-back square such
that $p$ is proper, $e$ is a closed embedding and $p^{-1}(X-e(A))\sr
X-e(A)$ is an isomorphism. 
\end{description}
\begin{remark}\rm
These cd-structures own their names to the fact that the behavior of
the functors of inverse image $f^*$ and $f^!$, which have upper
indexes, with respect to etale morphisms is very similar to the
behavior of the the functors of direct image $f_*$ and $f_!$, which
have lower indexes, with respect to proper morphisms.
\end{remark}
The topology associated with the upper cd-structure is called the
upper cd-topology. We will show below (see Proposition \ref{isnis})
that it coincides with the Nisnevich topology. In particular, an etale
morphism $f:X\sr Y$ is an upper covering if and only if for any $y$ in
$Y$ the fiber $p^{-1}(y)$ contains a $k_y$-rational point.  The
topology associated with the lower cd-structure is called the lower
cd-topology or proper cdh-topology. By Proposition \ref{lowerchar} a
proper morphism of schemes $p:X\sr Y$ is a lower cd-covering if and
only if for any point $y$ in $Y$ the fiber $p^{-1}(y)$ contains a
$k_y$-rational point.

The intersection of the upper and lower cd-structures is equivalent to
the {\em additive cd-structure} where a square is distinguished if it
is of the form
\begin{equation}
\label{addsq}
\begin{CD}
\emptyset@>>>Y\\
@VVV @VVV\\
X @>>> X\coprod Y
\end{CD}
\end{equation}
A presheaf $F$ is a sheaf in the topology associated with the additive
cd-structure if and only if
$$F(X\coprod Y)=F(X)\times F(Y)$$
and $F(\emptyset)=pt$. 

The union of the upper and lower cd-structures is the {\em combined
cd-structure}. A square is distinguished in it if it is an upper
distinguished or a lower distinguished square. Proposition \ref{isnis} and
the definition given in \cite[\S 4.1]{SusVoe2} imply that the
associated topology is the cdh-topology.

If we consider only squares where both $e$ and $p$ are monomorphisms
the upper and lower cd-structures become:
\begin{description}
\item[Plain upper cd-structure] or Zariski cd-structure where a square
of the form (\ref{eq1old}) is distinguished if both $p$ and $e$ are open
embeddings and $X=p(Y)\cup e(A)$. The associated topology is the
Zariski topology.
\item[Plain lower cd-structure] where a square of the form (\ref{eq1old})
is distinguished if both $p$ and $e$ are closed embeddings and
$X=p(Y)\cup e(A)$. The associated topology is the closed analog of the
Zariski topology. 
\end{description}
Any combination of the additive, upper, lower, plain upper and plain lower
cd-structures is called a standard cd-structure. There are nine
standard cd-structures: the five generating ones, the combined
cd-structure and the combinations of the plain upper with plain lower,
plain upper with lower and plain lower with upper cd-structures. They
form the following lattice where arrows indicate inclusions
$$
\begin{CD}
add @>>> p.up @>>> up\\
@VVV @VVV @VVV\\
p.low @>>> p @>>> p.low + up\\
@VVV @VVV @VVV\\
low @>>> low+ p.up @>>> cdh
\end{CD}
$$
The topology associated with the combination of the lower
and the plain upper cd-structures $low+p.up$ is considered in
\cite{GL}. 
\begin{lemma}
\llabel{stcompl} The standard cd-structures are complete on the
category of schemes or schemes of finite type over a base. In addition
the upper and plain upper cd-structures are complete on the category
of smooth schemes over a base and the lower and plain lower
cd-structures are complete on the category of proper schemes over a
base.
\end{lemma}
\begin{proof}
It follows immediately from \cite[Lemma 2.5]{HH1}.
\end{proof}
Let us show now that the standard cd-structures considered on the
category of schemes of finite dimension are bounded. A sequence of
points $x_0,\dots,x_d$ of a topological space $X$ is called an
increasing sequence (of length $d$) 
if $x_i\ne x_{i+1}$ and $x_i\in cl(\{x_{i+1}\})$ where
$cl(\{x_{i+1}\})$ is the closure of the point $x_{i+1}$ in $X$.
For a scheme $X$ define $D_d(X)$ as the class of open embeddings $j:U\sr X$ such that
for any $z\in X-U$ there exists an increasing sequence
$z=x_0,x_1,\dots,x_d$ of length $d$. The density structure defined by
the classes $D_d$ is called the {\em
standard density structure} on the category of schemes. It is locally
of finite dimension on the category of schemes of finite dimension and
the dimension of a scheme with respect to it is the dimension of the
corresponding topological space. 
\begin{lemma}
\llabel{conv3}
If $U,V\in D_{d}(X)$ then $U\cap V\in D_d(X)$.
\end{lemma}
\begin{lemma}
\llabel{conv2}
Let $U\in D_d(X)$ and $V$ be an open subscheme of $X$. Then $U\cap
V\in D_d(V)$.
\end{lemma}
\begin{proof}
Let $x$ be a point of $V$ outside of $U\cap V$. Considered as a point
of $X$ it has an increasing sequence $x=x_0,\dots,x_d$ with $x_i\in
X$. But since $x_0\in V$ we have $x_i\in V$ because $x_0\in cl(x_i)$
and $V$ is open.
\end{proof}
\begin{lemma}
\llabel{conv00}
Let $x_0,x_1,x_2$ be an increasing sequence on a 
scheme $X$ and $Z$ be a closed subset of $X$ such that $x_2$ lies
outside $Z$. Then there exists an increasing sequence $x_0,x_1',x_2$
such that $x_1'$ lies outside $Z$.
\end{lemma}
\begin{proof}
Replacing $X$ by the local scheme of $x_0$ in the closure of $x_2$ we
may assume 
that any point of $X$ contains $x_0$ in its closure and in turn lies
in the closure of $x_2$. It remains to show that the complement to $Z$
contains at least one point which is not equal to $x_2$. If it were
false we would have $x_2=X-Z$ i.e. $x_2$ would be a locally closed
point. This contradicts our assumption since by
\cite[5.1.10(ii)]{EGA4} a locally closed point on a locally Noetherian
scheme has dimension $\le 1$.
\end{proof}
\begin{lemma}
\llabel{conv0}
Let $X$ be a scheme, $U$ a dense open subset of $X$ and
$x_0,\dots,x_d$ any increasing sequence in $X$. Then there exists an
increasing sequence $x_0,x_1',\dots,x_d'$ such that $x_i'\in U$ for
all $i\ge 1$.
\end{lemma}
\begin{proof}
We may assume that $d>0$. If $x_d$ is contained in $U$ set
$x'_d=x_d$. Otherwise let $x'_d$ be a point of $U$ such that
$x_{d-1}\in cl(x_d)$ which exists since $U$ is dense. Since $x_d$ is
not in $U$,  $x_{d-1}$ is not in $U$ and thus $x'_d\ne x_{d-1}$ and
$x_0,x_1,\dots, x'_d$ is again an increasing sequence. Assume by
induction that we constructed $x_{i+1}',\dots,x_d'\in U$ such that
$x_0,\dots, x_i,x'_{i+1},\dots,x_d'$ is an increasing sequence. By
Lemma \ref{conv00} for any increasing sequence $y_0,y_1,y_2$ and a closed
subset $Z$ which does not contain $y_2$ there exists an increasing
sequence $y_0,y_1',y_2$ such that $Z$ does not contain $y_1$. Applying
this result to the sequence $x_{i-1},x_i,x'_{i+1}$ and 
$Z=X-U$ we construct $x_i'$.
\end{proof}
\begin{lemma}
\llabel{convv}
Let $X$ be a scheme and $Y$ a constructible subset in
$X$. Then any point $y'$ of the closure $cl(Y)$ of $Y$ in $X$ belongs to
the closure of a point $y$ of $Y$.
\end{lemma}
\begin{proof}
Since $Y$ is constructible it is of the form $Y=\cup_{i=1}^n Y_i$ where
each $Y_i$ is open in a closed subset of $X$ (see
e.g. \cite[Prop. 2.3.3]{EGA1}). It is clearly sufficient to prove our
statement for each $Y_i$. As a topological space $Y_i$ corresponds to
a Noetherian scheme. Thus there exists finitely many points $y_i'$ in
$Y$ such that any point of $Y$ is in the closure of one of the
$y_i$'s. If a point $y'$ in $cl(Y)$ has an open neighborhood $U$ which
does not contain any of the points $y_i$ then $U$ does not contain any
point of $Y$ which contradicts the assumption that $y\in cl(Y)$. Thus
$y$ belongs to the closure of $\{y_i\}$ which coincides with the union
of closures of points $y_i$ since there is finitely many of them.
\end{proof}
\begin{lemma}
\llabel{conv} Let $f:X\sr Y$ be a morphism of finite type of
Noetherian schemes and assume that there exists an open subset $U$ in
$Y$ such that $f^{-1}(U)$ is dense in $X$ and $f^{-1}(U)\sr U$ has
fibers of dimension zero. Then for any $d\ge 0$ and $V\in D_d(X)$
there exists $W\in D_{d}(Y)$ such that $f^{-1}(W)\subset V$.
\end{lemma}
\begin{proof}
We may clearly assume that $d>0$.  Let $Z=X-V$. We have to show that
$Y-cl(f(Z))\in D_{d}(Y)$ i.e. that for any $y$ in $cl(f(Z))$ there
exists an increasing sequence $y=y_0,\dots,y_d$ in $Y$. Since $f$ is
of finite type $f(Z)$ is constructible and in particular any point of
$cl(f(Z))$ is in the closure of a point in $f(Z)$ by Lemma
\ref{convv}. Thus we may assume that $y$ belongs to $f(Z)$
i.e. $y=f(x)$ where $x$ is in $Z$. By Lemma \ref{conv0} we can find an
increasing sequence $x=x_0,x_1,\dots, x_d$ for $x$ such that for $i>0$
we have $x_i\in f^{-1}(U)$. Then $y=f(x_0),\dots,f(x_d)$ is an
increasing sequence i.e. $f(x_i)\ne f(x_{i+1})$. Indeed for $i>0$ it
follows from the fact that the fibers of $f$ over $U$ are of dimension
zero. For $i=0$ we have two cases. If $f(x_0)\in U$ then the same
argument as for $i>0$ applies. If $f(x_0)$ is not in $U$ then
$f(x_0)\ne f(x_1)$ since $f(x_1)\in U$.
\end{proof}
\begin{proposition}
\llabel{upred}
The upper cd-structure and the plain upper cd-structures on the
category of Noetherian schemes of finite dimension are bounded with
respect to the standard density structure.
\end{proposition}
\begin{proof}
We will only consider the upper cd-structure. The plain case is
similar. Let us show that any upper distinguished square is reducing
with respect to the standard density structure (see \cite[Definition
2.19]{HH1}).  Let our square be of the form
\begin{equation}
\begin{CD}
W @>j_V>> V\\
@VVV @VVpV\\
U @>j>> X
\end{CD}
\end{equation}
and let $W_0\in D_{d-1}(W)$, $U_0\in D_{d}(U)$, $V_0\in
D_d(V)$. Applying Lemma \ref{conv} to the morphism $j\coprod p$ we can
find $X_0\in D_d(X)$ such that $j(U_0)\cup p(V_0)\subset
X_0$. Replacing $X$ with $X_0$ and applying Lemma \ref{conv2} we may
assume that $U_0=U$ and $V_0=V$. Let $Z=W-W_0$, $C=X-U$ and set
$X'=X-(C\cap cl(pj_V(Z)))$.  Let us show that the
square
\begin{equation}
\begin{CD}
W_0@>>>j_V(W_0)\\
@VVV @VVV\\
U@>>>X'
\end{CD}
\end{equation}
is upper distinguished. It is clearly a pull back square, the right
vertical arrow is etale and the lower horizontal one is an open
embedding. It is also obvious that $p^{-1}(X'-U)\cap j_V(W_0)=(X'-U)$. 
To finish the proof it remains to show that $X'\in D_d(X)$. Let $x$ be
a point of $X$ outside of $X'$ i.e. a point of $C\cap
cl(pj_V(Z))$. Since $pj_V(Z)\cap C=\emptyset$ there exists $x'=pj_V(y)\in
pj_V(Z)$ such that $x\in cl(x')$ and $x'\ne x$. Let
$y=y_0,\dots,y_{d-1}$ be an increasing sequence for $y$ in $W$ which
exists since $W_0\in D_{d-1}(W)$. The morphism $q=pj_V$ has fibers of
dimension zero and therefore $q(y_0),\dots,q(y_{d-1})$ is an
increasing sequence for $x'$. Thus we get an increasing sequence
$x,q(y_0),\dots,q(y_{d-1})$ for $x$ of length $d$.
\end{proof}
\begin{proposition}
\llabel{lowerred}
The lower cd-structure and the plain lower cd-structures on the
category of Noetherian schemes of finite dimension are bounded with
respect to the standard density structure.
\end{proposition}
\begin{proof}
We will only consider the case of the lower cd-structure. The plain
case is similar. Consider a lower distinguished square
\begin{equation}
\llabel{sqag}
Q=\left(
\begin{CD}
B @>i_Y>> Y\\
@VVV @VVpV\\
A @>i>> X
\end{CD}
\right)
\end{equation}
If we replace $Y$ by the scheme-theoretic closure of the open
subscheme $p^{-1}(X-A)$ we get another lower distinguished square
which is a refiniment of the original one. This square satisfies the
condition of Lemma \ref{l8.4.1} and therefore it is reducing.
\end{proof}
\begin{lemma}
\llabel{l8.4.1} A lower distinguished square of the form (\ref{sqag})
such that the subset $p^{-1}(X-A)$ is dense in $Y$ is reducing with
respect to the lower cd-structure.
\end{lemma}
\begin{proof}
Let $Y_0\in D_{d}(Y)$, $A_0\in D_d(A)$, $B_0\in
D_{d-1}(B)$. Applying Lemma \ref{conv} to $p$ and $U=X-A$ we conclude
that there exists $X_0\in D_{d}(X)$ such that $p(Y_0)\subset
X_0$. Applying the same lemma to $i$ we find an open subset $X_1\in
D_d(X)$ such that $i(A_0)\subset X_1$. Then by Lemma \ref{conv3}
$X_1\cap X_0\in D_d(X)$ and replacing $X$ by $X_1\cap X_0$ and using
Lemma \ref{conv2} we may assume that $A_0=A$ and $Y_0=Y$. Let
$X'=X-pi_Y(B-B_0)$. To finish the proof it is enough to check that
$X'\in D_d(X)$ and define $Q'$ as the pull-back of $Q$ to
$X'$. According to Lemma \ref{conv} applied again to $p$ and $U=X-A$
it is enough to check that $Y-i_Y(B-B_0)\in D_{d}(Y)$. Since $B_0\in
D_{d-1}(B)$ and $i_Y$ is a closed embedding it is enough to check that
$Y-i_B(B)$ is dense in $Y$. This follows from our  assumption since
$Y-i_B(B)=p^{-1}(X-A)$. 
\end{proof}
Since all generating cd-structures on the category of Noetherian
schemes are bounded by the {\em same} density structure any their
combination is also bounded by the same density structure. We get the
following result.
\begin{proposition}
\llabel{cdhis}
The standard cd-structures on the category of Noetherian schemes
of finite dimension are bounded.
\end{proposition}
Finally let us show that all the standard cd-structures are
regular. It is clearly sufficient to consider the ``generating''
cd-structures. Then any combination of them will also be regular. 
\begin{lemma}
\llabel{uploreg}
The additive, upper, plain upper, lower and plain lower cd-structures are
regular. 
\end{lemma}
\begin{proof}
The additive case is obvious.  Let us show that the upper, palin
upper, lower and plain lower cd-structures satisfy the conditions of
\cite[Lemma 2.11]{HH1}. The first two conditions are obvious. Consider
the third condition in the upper case. The square
\begin{equation}
\llabel{eq2}
d(Q)=\left(
\begin{CD}
B@>>>Y\\
@VVV @VVV\\
B\times_A B@>>>Y\times_X Y
\end{CD}
\right)
\end{equation}
is a pull-back square. Since $p$ is etale, and in particular
unramified, the diagonal $Y\sr Y\times_X Y$ is an open embedding. The
morphism $B\times_A B\sr Y\times_X Y$ is an open embedding because $e$
is an open embedding. The condition that $p^{-1}(X-e(A))\sr
X-e(A)$ is a universal homeomorphism implies that for a pair of
geometric points $y_1,y_2$ of $Y$ such that $p(y_1)=p(y_2)\in X-e(A)$
one has $y_1=y_2$. Therefore, 
$$Y\times_X Y= (B\times_A B)\cup Y$$
i.e. (\ref{eq2}) is a (plain) upper distinguished square.

Consider the third condition in the lower case. The square (\ref{eq2})
is a pull-back square. Since $p$ is proper, and in particular
separated, the diagonal $Y\sr Y\times_X Y$ is a closed embedding. The
morphism $B\times_A B\sr Y\times_X Y$ is a closed embedding because
$e$ is a closed embedding. The condition that
$p^{-1}(X-e(A))\sr X-e(A)$ is a universal homeomorphism implies that
for a pair of geometric points $y_1,y_2$ of $Y$ such that
$p(y_1)=p(y_2)\in X-e(A)$ one has $y_1=y_2$. Therefore,
$$Y\times_X Y= (B\times_A B)\cup Y$$
i.e. (\ref{eq2}) is a (plain) lower distinguished square.
\end{proof}
\begin{definition}
\llabel{uppersplit}
Let $f:\tilde{X}\sr X$ be a morphism of schemes. A splitting
sequence for $f$ is a sequence of closed embeddings 
$$\emptyset=Z_{n+1}\sr Z_n\sr\dots\sr Z_1\sr Z_0=X$$
such that for any $i=0,\dots,n$ the projection
$$(Z_i-Z_{i+1})\times_X \tilde{X}\sr (Z_i-Z_{i+1})$$ has a section.
\end{definition}
\begin{lemma}
\llabel{l8.5.1} A morphism of finite type of Noetherian schemes
$f:\tilde{X}\sr X$ has a splitting sequence if and only if for any
point $x$ of $X$ there exists a point $\tilde{x}$ of $\tilde{X}$ such
that $f(\tilde{x})=x$ and the corresponding morphism of the residue
fields is an isomorphism.
\end{lemma}
\begin{proof}
The ``only if'' part is obvious. The ``if'' part follows easily by
the Noetherian induction (cf. \cite[Lemma 3.1.5]{MoVo}).
\end{proof}
\begin{proposition}
\llabel{isnis} An etale morphism $f:\tilde{X}\sr X$ is a covering in
the upper cd-topology if and only if for any point $x$ of $X$ there
exists a point $\tilde{x}$ of $\tilde{X}$ such that $f(\tilde{x})=x$
and the corresponding morphism of the residue fields is an
isomorphism.
\end{proposition}
\begin{proof}
Since the upper cd-structure is complete any upper cd-covering has a
refinement which is a simple covering which immediately implies the
``only if'' part of the proposition. To prove the ``if'' part we have
to show, in view of Lemma \ref{l8.5.1}, that any etale morphism
$f:\tilde{X}\sr X$ which has a splitting sequence $Z_n\sr\dots\sr
Z_0=X$ is an upper cd-covering.  We will construct an upper
distinguished square of the form (\ref{eq1}) based on $X$ such that
the pull-back of $f$ to $Y$ has a section and the pull-back of $f$ to
$A$ has a splitting sequence of length less than $n$. The result then
follows by induction on $n$. We take $A=X-Z_n$. To define $Y$ consider
the section $s$ of $f_n:\tilde{X}\times_X Z_n\sr Z_n$ which exists by
definition of a splitting sequence. Since $f$ is etale and in
particular unramified the image of $s$ is an open subscheme. Let $W$
be its complement. The morphism $\tilde{X}\times_X Z_n\sr \tilde{X}$
is a closed embedding thus the image of $W$ is closed in
$\tilde{X}$. We take $Y=\tilde{X}-W$. One verifies immediately that
the pull-back square defined by $A\sr X$ and $Y\sr X$ is upper
distinguished. The pull-back of $f$ to $Y$ has a section and the
pull-back of $f$ to $A$ has a splitting sequence of length $n-1$. This
finishes the proof of the proposition.
\end{proof}
Proposition \ref{isnis} implies that the topology associated with the
upper cd-structure on the category of Noetherian schemes is the
Nisnevich topology.
\begin{proposition}
\llabel{lowerchar} 
A proper morphism
$f:\tilde{X}\sr X$ is a covering in the lower cd-topology if and
only if for any point $x$ of $x$ there exists a point $\tilde{x}$ of
$\tilde{X}$ such that $f(\tilde{x})=x$ and the corresponding morphism
of the residue fields is an isomorphism.
\end{proposition}
\begin{proof}
Since the lower cd-structure is complete any lower cd-covering has a
refinement which is a simple covering which immediately implies the
``only if'' part of the proposition. To prove the ``if'' part we have
to show, in view of Lemma \ref{l8.5.1}, that any proper morphism
$f:\tilde{X}\sr X$ which has a splitting sequence $Z_n\sr\dots\sr
Z_0=X$ is a lower cd-covering. We will construct a lower distinguished
square of the form (\ref{eq1}) based on $X$ such that the pull-back of
$f$ to $Y$ has a section and the pull-back of $f$ to $A$ has a
splitting sequence of length less than $n$. The result then follows by
induction on $n$. We take $A=Z_1$. To define $Y$ consider the section
$s$ of $f_n:\tilde{X}\times_X (X-Z_1)\sr (X-Z_1)$ which exists by
definition of a splitting sequence. Since $f$ is proper and in
particular separated, the image of $s$ is a closed subscheme. Let $W$
be its complement. The morphism $\tilde{X}\times_X (X-Z_1)\sr
\tilde{X}$ is an open embedding thus the image of $W$ is open in
$\tilde{X}$. We take $Y=\tilde{X}-W$. One verifies immediately that
the pull-back square defined by $A\sr X$ and $Y\sr X$ is lower
distinguished. The pull-back of $f$ to $Y$ has a section and the
pull-back of $f$ to $A$ has a splitting sequence of length $n-1$. This
finishes the proof of the proposition.
\end{proof}

\subsection{Motivic homotopy categories}
Recall that in \cite{MoVo} we defined for any site $T$ with an
interval $I$ a category $H(T,I)$ which we called the homotopy category
of $(T,I)$. Applying this definition to a category of schemes with
some standard topology and taking $I$ to be the affine line one
obtains different motivic homotopy categories. Among these homotopy
categories the one denoted in \cite{MoVo} by $H(S)$ and corresponding
to the category of smooth schemes over $S$ with the Nisnevich or upper
cd-topology seems to play a distinguished role. In this section we prove a
number of results which provide a new description for the motivic
homotopy categories in the standard topologies and in particular for
the category $H(S)$. We start with some results applicable to all
sites with interval with good enough completely decomposable
topologies.

Let $C$ be a category with a complete regular bounded cd-structure $P$
(see \cite{HH1}) and an interval $I$ (see \cite[\S 2.3]{MoVo}). Assume
in addition that $C$ has a final object and that for any $X$ in $C$
the product $X\times I$ exists.  We can form the homotopy category of
$(C,P,I)$ in two ways. First, we may define a new cd-structure $(P,I)$
whose distinguished squares are the distinguished squares of $C$ and
squares of the form
$$
\begin{CD}
\emptyset@>>>\emptyset\\
@VVV @VVV\\
X\times I @>>> X
\end{CD}
$$
where $X$ runs through all objects of $C$ and consider the homotopy
category $H(C,P,I)$ of this cd-structure i.e. the localization of
$\Delta^{op}PreShv(C)$ with respect to $cl_{\bdl}(W_{(P,I)}\cup
W_{proj})$ where $cl_{\bdl}(-)$ is the $\bdl$-closure defined in
\cite{HH0}. On the other hand we may consider the homotopy category ${
H}(C_{t_P},I)$ of the site with interval $(C_{t_P},I)$ as defined in
\cite{MoVo}. We are going to show that if $P$ is complete, regular and
bounded these two constructions agree. Since the comparison theorem of
the next section is formulated in terms of pointed categories we use
pointed context to formulate the results of this section as well. One
can easily see that the same arguments can be used to prove the
corresponding results in the free context. Recall, that for a functor
$\Phi$ we denote by $iso(\Phi)$ the class of morphisms $f$ such that
$\Phi(f)$ is an isomorphism.
\begin{proposition}
\llabel{swi} Let $P$ be a complete, regular and bounded
cd-structure. Then the functor
\begin{equation}
\llabel{isl}
\Phi:\Delta^{op}PreShv_{\bullet}(C)\sr {H}_{\bullet}(C_{t_P},I)
\end{equation}
is a localization and 
$$iso(\Phi)=cl_{\bdl}((W_{(P,I)})_+\cup
W_{proj})=cl_{\bdl}((W_{P})_+\cup (W_I)_+\cup
W_{proj})$$
where $W_P$ is the class of generating weak equivalences of $P$, $W_I$
the class of all projections $X\times I\sr X$ for $X\in C$ and
$W_{proj}$ is the class of projective weak equivalences of (pointed)
simplicial presheaves.
\end{proposition}
\begin{proof}
Recall that the category ${H}_{\bullet}(C_{t_P},I)$ is defined as the
localization of the category of pointed simplicial sheaves on $C$ in
the $t_P$-topology with respect to $I$-weak equivalences (see
\cite[Def. 2.3.1]{MoVo}). Since the category of pointed simplicial
sheaves is a localization of the category of pointed simplicial
presheaves the functor $\Phi$ is a localization. It remains to check
that a morphism of pointed simplicial presheaves $f$ belongs to
$cl_{\bdl}((W_{P}\cup W_I)_+\cup W_{proj})$ if and only if the
associated morphism of sheaves is an I-weak equivalence.  The
morphisms associated with elements of $(W_I)_+$ are I-weak
equivalences by definition and since our cd-structure is regular so
are the morphisms associated with elements of $(W_P)_+$. The ``only
if'' part follows now from Lemma \ref{isdlta}. Let $N$ be the set of
morphisms of the form $(\emptyset\sr U)_+$ for $U\in C$ and let
$Ex=Ex_{W_P\cup W_I, N}$ be the functor constructed in
\cite[Proposition 2.2.12]{HH0} such that for any $X$ the morphism
$X\sr Ex(X)$ is in $cl_{\bdl}((W_P\cup W_I)_+)$, the simplicial sets
$Ex(X)(U)$ are Kan and for any $f:Y\sr Y'$ in $(W_P\cup W_I)_+$ the
map $S(Y',Ex(X))\sr S(Y, Ex(X))$ defined by $f$ is a weak
equivalence. In view of Lemma \ref{isdlta} the morphisms associated
with the morphisms $X\sr Ex(X)$ are $I$-weak equivalences. To finish
the proof it is sufficient to show that any morphism $f:Ex(X)\sr
Ex(Y)$ such that the associated morphism of sheaves is an $I$-weak
equivalence is a projective weak equivalence. Let $Ex_{JJ}(X)$ be a
fibrant replacement of $a(Ex(X))$ in the Jardine-Joyal closed model
structure on the category of sheaves in $t_P$. The morphism $Ex(X)\sr
Ex_{JJ}(X)$ is a local weak equivalence with respect to $t_P$ and,
since both objects are flasque with respect to $P$, \cite[Lemma
3.5]{HH1} implies that it is a projective weak equivalence. Together
with the fact that $Ex(X)(U)\sr Ex(X)(U\times I)$ is a weak
equivalence for any $U$ in $C$ this implies that $Ex_{JJ}(X)$ is
$I$-local. If $f:Ex(X)\sr Ex(Y)$ is a morphism such that $a(f)$ is an
$I$-weak equivalence then the corresponding morphism $Ex_{JJ}(X)\sr
Ex_{JJ}(Y)$ is an $I$-weak equivalence and, therefore, a local weak
equivalence with respect to $t_P$. We conclude that $f$ is a local weak
equivalence and, using again the fact that $Ex(-)$ are flasque and
\cite[Lemma 3.5]{HH1}, we conclude that $f$ is a projective weak
equivalence.
\end{proof}
\begin{lemma}
\llabel{isdlta} The class of morphisms $f$ such that $a(f)$ is a
pointed I-weak equivalence is $\bdl$-closed.
\end{lemma}
\begin{proof}
The associated sheaf functor commutes with coproducts and therefore it
is enough to show that the class of (pointed) I-weak equivalences in
$\Delta^{op}Shv_{t_P, \bullet}(C)$ is $\bdl$-closed.  The fact that
the class of I-weak equivalences is closed under coproducts follows
from its definition. It also follows from its definition that this
class satisfies the first two conditions of the definition of a
$\Delta$-closed class (see \cite{HH0}). To verify the third condition
observe that since homotopy colimits and diagonals are defined on
simplicial presheaves by applying the corresponding construction for
simplicial sets over each object of the category, \cite[Ch.XII
4.3]{BKan} implies that for a bisimplicial sheaf $B$ there is a
natural weak equivalence $hocolim B_i\sr \Delta B$ where $B_i$ are the
rows (or columns) of $B$. The condition follows now from \cite[Lemma
2.2.12]{MoVo}. The fact that the class of I-weak equivalences is
closed under colimits of sequences of morphisms is proved in
\cite[Cor. 2.2.13(2)]{MoVo}.
\end{proof}
\begin{cor}
\llabel{c8.3.1}
Under the assumptions of Proposition \ref{swi} the functor 
$$\Phi:\Delta^{op}C^{\coprod}_+\sr H(C,P,I)$$
is a localization and $iso(\Phi)=cl_{\bdl}((W_P\cup W_I)_+)$. 
\end{cor}
\begin{proof}
This is a particular case of \cite[Corollary 4.3.8]{HH0}.
\end{proof}
\begin{cor}
\llabel{stricteq} Let $C$ be any category with an interval and $f:X\sr
Y$ be a strict $I$-homotopy equivalence in
$\Delta^{op}C^{\coprod}_+$. Then $f$ belongs to $cl_{\bdl}((W_I)_+)$.
\end{cor}
\begin{proof}
By \cite[Lemma 2.3.6]{MoVo} applied to the site $C$ with the trivial
topology we know that any strict homotopy equivalence is an $I$-weak
equivalence. Applying Corollary \ref{c8.3.1} to the case of the
empty cd-structure on the category obtained from $C$ by the addition
of an initial object we conclude that any strict $I$-homotopy
equivalence belongs to $cl_{\bdl}((W_I)_+)$
\end{proof}
\begin{remark}\rm
Results of \cite{HH0} imply that the category $H_{\bullet}(C,P,I)$ is
the homotopy category of a simplicial almost finitely generated closed
model structure on the category of simplicial presheaves on $C$ where
cofibrations are the projective cofibrations. If $P$ is complete,
bounded and regular we can also realize it as the homotopy category of
an almost finitely generated closed model structure on the category of
simplicial sheaves on $C$ in the $t_P$-topology. Indeed, Proposition
\ref{swi} states that, under these assumptions, the category
$H_{\bullet}(C,P,I)$ is equivalent to the localization of the
homotopy category of pointed simplicial sheaves in $t_P$ with respect
to the class of $(W_I)_+$-local equivalences (see \cite{Hirs}). The
homotopy category of simplicial sheaves is the homotopy category of
the Brown-Gersten closed model structure which is cellular by
\cite[Proposition 4.7]{HH1}. Therefore, $H_{\bullet}(C,P,I)$ is
equivalent to the homotopy category of the left Bousfield localization
of the Brown-Gersten closed model structure with respect to $(W_P)_+$
which exists by \cite{Hirs} and one verifies easily that it is almost
finitely generated.
\end{remark}
Specializing these general theorems to the case of the motivic
homotopy categories and using the properties of the standard
cd-structures proved in the first section we get the following
results.
\begin{proposition}
\llabel{v1} Let $P$ be a standard cd-structure on the category
$Sch/S$. Then the functor
\begin{equation}
\llabel{v1eq}
\Delta^{op}(Sch/S)^{\coprod}_+\sr H_{\bullet}((Sch/S)_{t_P},\af) 
\end{equation}
is the
localization with respect to the smallest $\bdl$-closed
class which contains morphisms of the form $(p_Q:K_Q\sr X)_+$ for 
distinguished squares $Q$ and morphisms of the form $(X\times\af\sr X)_+$
for schemes $X$.
\end{proposition}
\begin{proposition}
\llabel{v2}
Let $P$ be a standard cd-structure which is contained in the upper
cd-structure. Then the functor
\begin{equation}
\llabel{v2eq}
\Delta^{op}(Sm/S)^{\coprod}_+\sr H_{\bullet}((Sm/S)_{t_P},\af)
\end{equation}
is the localization with respect to the smallest $\bdl$-closed class
which contains morphisms of the form $(p_Q:K_Q\sr X)_+$ for upper
distinguished squares $Q$ and morphisms of the form $(X\times\af\sr
X)_+$ for smooth schemes $X$.
\end{proposition}
In many cases the category of all schemes (resp. all smooth schemes)
on the left hand side of (\ref{v1eq}) and (\ref{v2eq}) can be replaced
by smaller subcategories. For either upper or lower cd-structure all
schemes can be replaced by quasi-projective schemes (for the lower
cd-structure one uses the Chow lemma to show that this is
allowed). For plain upper or stronger cd-structure smooth schemes can
be replaced by smooth quasi-affine schemes etc.

\subsection{The comparison theorem}
Let $k$ be a field. We have an obvious functor of pointed motivic
homotopy categories 
\begin{equation}
\llabel{niscdh}
{ H}_{\bullet}((Sm/k)_{Nis},\af)\sr {
H}_{\bullet}((Sch/k)_{cdh},\af)
\end{equation}
which we denote by ${\bf L}\pi^*$ because it is the inverse image
functor defined by the continuous map of sites
$$\pi:(Sch/k)_{cdh}\sr (Sm/k)_{Nis}$$
For a morphism $f$ in the pointed homotopy category we denote by
$$\Sigma^1_s(f)=f\wedge Id_{S^1_s}$$
the first simplicial suspension of $f$.
Let us recall the following definition given in \cite{FV}.
\begin{definition}
\llabel{res}
A field $k$ is said to admit resolution of singularities if the
following two conditions hold:
\begin{enumerate}
\item for any reduced scheme of finite type $X$ over $k$ there exists
a proper morphism $f:\tilde{X}\sr X$ such that $\tilde{X}$ is smooth
and $f$ has a section over a dense open subset of $X$
\item for any smooth scheme $X$ over $k$ and
a proper surjective morphism $Y\sr X$ which has a section over a dense
open subset of $X$ there exists a sequence of blow-ups with smooth
centers $X_n\sr X_{n-1}\sr\dots\sr X_0=X$ and a morphism $X_n\sr Y$
over $X$.
\end{enumerate}
\end{definition}
Note that any field satisfying the conditions of Definition \ref{res}
is perfect. 
\begin{theorem}
\llabel{m124} Let $k$ be a field which admits resolution of
singularities. Then the functor ${\bf L}\pi^*$ is a localization and
for any $f$ in $iso({\bf L}\pi^*)$ the morphism $\Sigma^1_s(f)$ is an
isomorphism.
\end{theorem}
\begin{proof}
Define the smooth blow-up cd-structure on the category $Sm/k$ of
smooth schemes over $k$ as the collection of pull-back squares of the
form (\ref{eq1}) such that $e$ is a closed embedding and $p$ is the
blow-up with the center in $e(A)$.
\begin{lemma}
\llabel{iscomplete} Let $k$ be a field which admits resolution
of singularities. Then the smooth blow-up cd-structure on the category
of smooth schemes over $k$ is complete.
\end{lemma}
\begin{proof}
To show that a cd-structure is complete it is sufficient to show that
for any distinguished square of the form (\ref{eq1}) and any morphism
$f:X'\sr X$ the sieve $f^*(e,p)$ contains the sieve generated by a
simple covering (see \cite[Lemma 2.4]{HH1}). Let us prove it by
induction on $dim(X')$. If $dim(X')=0$ the sieve $f^*(e,p)$ contains
an isomorphism. Assume that the statement is proved for $dim(X')<d$
and let $X'$ be of dimension $d$. The map $X'\times_X (A\coprod Y)\sr
X'$ is proper and has a section over a dense open subset of $X'$. Thus
by the resolution of singularities assumption we have a sequence of
blow-ups with smooth centers $X'_n\stackrel{p_{n-1}}{\sr}
X'_{n-1}\stackrel{p_{n-2}}{\sr} \dots\stackrel{p_{0}}{\sr}X'_0=X'$
such that the pull-back of $(e,p)$ to $X'_n$ contains an isomorphism
and in particular a sieve generated by a simple covering. Assume by
induction that the pull-back of $(e,p)$ to $X_{i}'$ contains a sieve
generated by a simple covering $\{r_j:U_{j}\sr X_{i}'\}$ and let us
show that the same is true for $X_{i-1}'$. Let $e_{i-1}:Z'_{i-1}\sr
X'_{i-1}$ be the center of the blow-up $X_{i}'\sr X_{i-1}'$. The
restriction of $(e,p)$ to $Z_{i-1}'$ contains a sieve generated by a
simple covering $\{s_l:V_{l}\sr Z_{i-1}'\}$ since
$dim(Z_{i-1}')<d$. Thus the restriction of $(e,p)$ to $X_{i-1}'$
contains the sieve generated by $\{p_{i-1}r_j, e_{i-1}s_l\}$ which is
a simple covering by definition.
\end{proof}
\begin{lemma}
\llabel{smred}
The smooth blow-up cd-structure on the category of smooth schemes over
any field is bounded with respect to the standard density structure.
\end{lemma}
\begin{proof}
The same arguments as in the proof of Lemma \ref{l8.4.1} show that any
distinguished square of the smooth blow-up cd-structure is reducing
with respect to the standard density structure.
\end{proof}
\begin{lemma}
\llabel{regbu}
The smooth blow-up  cd-structure on the category of smooth schemes over
any field is regular.
\end{lemma}
\begin{proof}
The first two conditions of \cite[Definition 2.10]{HH1} are obviously
satisfied. To prove the third one we have to show that for a
distinguished square of the form (\ref{eq1}) the map
of representable sheaves of the form 
\begin{equation}
\llabel{eq8.5.1}
\rho(Y)\coprod \rho(B)\times_{\rho(A)}\rho(B)\sr
\rho(Y)\times_{\rho(X)}\rho(Y)
\end{equation}
is surjective. Since any smooth scheme has a covering in our topology
by connected smooth schemes it is sufficient to show that the map of
presheaves corresponding to (\ref{eq8.5.1}) is surjective on sections
on smooth connected schemes. Let $U$ be a smooth connected scheme and
$f,g:U\sr Y$ be a pair of morphisms such that $p\circ f=p\circ g$. The
scheme $Y\times_X Y$ is the union of two closed subschemes namely the
diagonal $Y$ and $B\times_A B$ (see the proof of the lower case in
Lemma \ref{uploreg}). Since $U$ is smooth and connected it is
irreducible and therefore the closure of the image of $f\times g$ in
$Y\times_X Y$ is irreducible. This implies that the image belongs to
either $Y$ or $B\times_A B$ and since $U$ is smooth and in particular
reduced the morphism $f\times_X g$ lifts to $Y$ or to $B\times_A
B$.
\end{proof}
Consider the topology $scdh$ associated with the sum of the smooth
blow-up cd-structure and the upper cd-structure on the category of
smooth schemes over $S$. Since the sum of two cd-structures bounded by
the same density structure is bounded, Proposition \ref{upred} and
Lemma \ref{smred} imply that this cd-structure is bounded by the
standard density structure on $Sm/k$. Since the sum of two regular
cd-structures is regular, Lemma \ref{uploreg} and Lemma \ref{regbu}
imply that it is regular. Since the sum of two complete cd-structures
is complete, Lemma \ref{iscomplete} and Lemma \ref{stcompl} imply that
if $k$ admits resolution of singularities then this cd-structure is
complete. Therefore, Corollary \ref{c8.3.1} implies that the functor
$$\Delta(Sm/k)^{\coprod}_{+}\sr {
H}_{\bullet}((Sm/k)_{scdh},\af)$$
is the localization with respect to the smallest
$\bdl$-closed class which contains morphisms of the form 
$(p_Q:K_Q\sr X)_+$, where $Q$ is an upper distinguished square or a
smooth blow-up square, and the projections $(X\times\af\sr X)_+$. On
the other hand the continuous map of sites
$$(Sch/k)_{cdh}\sr (Sm/k)_{scdh}$$
defined by the inclusion of categories $Sm/k\sr Sch/k$ defines the
inverse image
functor 
\begin{equation}
\llabel{invim}
{
H}_{\bullet}((Sm/k)_{scdh},\af)\sr {
H}_{\bullet}((Sch/k)_{cdh},\af)
\end{equation}
and we have a commutative diagram
$$
\begin{CD}
\Delta(Sm/k)^{\coprod}_{+}@>>>H_{\bullet}((Sm/k)_{Nis},\af)\\
@VVV @VV{\bf L}\pi^*V\\
{H}_{\bullet}((Sm/k)_{scdh},\af)@>>>{
H}_{\bullet}((Sch/k)_{cdh},\af)
\end{CD}
$$
\begin{lemma}
\llabel{easy}
If $k$ admits resolution of singularities the inverse image functor 
$Shv_{scdh}(Sm/k)\sr Shv_{cdh}(Sch/k)$ is an equivalence.
\end{lemma}
\begin{proof}
The resolution of singularities assumption implies that any object of
$Sch/k$ has a cdh-covering by objects of $Sm/k$ and that any
cdh-covering of an object of $Sm/k$ has a refinement which is a
scdh-covering. These to facts together imply that the inverse and the
direct image functors define equivalences of the corresponding
categories of sheaves (see \cite{}).
\end{proof}
Lemma \ref{easy} implies that the functor (\ref{invim}) is an
equivalence. Thus we conclude that the functor ${\bf L}\pi^*$ a
localization. By \cite[Lemma 3.4.13]{HH0} any morphism in
$H_{\bullet}((Sm/k)_{Nis},\af)$ is isomorphic to the image of a
morphism in $\Delta(Sm/k)^{\coprod}_{+}$ which implies that
$$iso({\bf L}\pi^*)=cl_{\bdl}(W_{scdh,+}\cup W_{\af,+})$$
For a class $E$ in $\Delta^{op}C^{\coprod}_{\bullet}$ and a pointed
simplicial set $K$ one has
$$cl_{\bdl}(E)\wedge Id_K\subset cl_{\bdl}(E\wedge
Id_K)$$ 
and, therefore,
$$\Sigma^1_s(iso({\bf L}\pi^*))\subset
cl_{\bdl}(\Sigma^1_s(W_{scdh,+})\cup \Sigma^1_s(W_{\af_+}))$$
Elements of $\Sigma^1_s(W_{\af_+})$ are $\af$-weak equivalences for
any topology, elements of $\Sigma^1_s(W_{scdh,+})$ are $\af$-weak
equivalences for the Nisnevich topology by \cite[Remark
3.2.30]{MoVo}. Together with Lemma \ref{isdlta} it implies that 
the class $\Sigma^1_s(iso({\bf L}\pi^*))$ consists of $\af$-weak
equivalences in the Nisnevich topology.
\end{proof}
\begin{cor}
\llabel{cortwo}
Let $k$ be as above and $X$ and $Y$ be pointed simplicial sheaves on
$(Sm/k)_{Nis}$ such that $Y$ is $\af$-weak equivalent to the
simplicial loop space of an $\af$-local object. Then the map
$$Hom(X,Y)\sr Hom({\bf L}\pi^*(X),{\bf L}\pi^*(Y)),$$ where the
morphisms on the left hand side are in
$H((Sm/k)_{Nis},\af)$ and on the right hand side in
$H_{\bullet}((Sch/k)_{cdh},\af)$, is bijective.
\end{cor}
\bibliography{alggeom}
\bibliographystyle{plain}
\end{document}